\documentclass[notitlepage, 12pt]{report}

\title{When is a Tensor Product of Circulant Graphs Circulant?}
\author{J.C.\ George\\Department of Mathematics\\Tougaloo College
\\Tougaloo, MS  39174 \and R.S.\ Sanders\\Department of Mathematics
\\Buffalo State College\\Buffalo, NY  14222}
\date{}
\newcommand{\Aut}{\mbox{Aut}}
\newcommand{\lcm}{\mbox{lcm}}
\newcommand{\qed}{\rule[-.5mm]{2mm}{3mm}}
\newtheorem{theorem}{Theorem}
\newtheorem{corollary}[theorem]{Corollary}
\newtheorem{lemma}[theorem]{Lemma}
\newenvironment{proof}[1][Proof:]
     {\begin{sloppypar} {\noindent\bf #1} }
     { \qed \end{sloppypar}} 
\begin{document}
\maketitle
\begin{abstract}
The tensor product, also called direct, categorical, or Kronecker 
product of graphs, is one of the least-understood graph products.  In 
this paper we determine partial answers to the question given in 
the title, thereby significantly extending results of Broere and 
Hattingh (see \cite{B&H}).  We characterize 
completely those pairs of complete graphs whose tensor products are 
circulant.  We establish that if the orders of these circulant graphs 
have greatest common divisor of 2, the product is circulant whenever both 
graphs are bipartite.  We also establish that it is possible for one of 
the two graphs not to be circulant and the product still to be circulant.
\end{abstract}

Throughout this paper, we will assume graphs to be connected unless 
otherwise stated, and without multiple edges.  However, we will permit 
graphs to have loops at some vertices.

The \emph{tensor product} $G \otimes H$ of two graphs $G$ and $H$ is 
the graph with vertex set $V(G) \times V(H)$ and edge set consisting of 
those pairs of vertices $(g,h), (g',h')$ where $g$ is adjacent to $g'$ 
and $h$ is adjacent to $h'$.

This product, also called the Kronecker product, weak product, direct 
product, categorical product, and conjunction,  has been studied for 
decades.  Basic properties of the tensor product may be found, for 
example, in \cite{PMW, H&W, Robin1, Robin2}.  Despite this study, many 
basic properties of the product are unknown or only partly understood.  
In this paper we determine when the tensor product of circulant graphs 
is circulant.  This research expands upon the work of Broere and 
Hattingh (\cite{B&H}).

A \emph{circulant graph} on $n$ vertices is a graph $C_nS$, where $S$ 
is a set of integers, with vertices labeled by integers from 0 to 
$n-1$.  Two vertices are adjacent in the graph if their labels differ 
by one of the values in $S$, mod $n$.  For convenience, we assume that 
all the entries in $S$ may be taken to be in the range 
$[0, \lfloor n/2 \rfloor -1]$.

We begin by stating a few elementary facts about circulant graphs
without 
proof.  Some of these facts may be found in \cite{B&H}, and others 
follow immediately from the definition.

(C1) A graph is circulant if and only if its automorphism group 
contains a transitive cyclic subgroup.

(C2) A graph is circulant if and only if its complement is circulant.

(C3) The graph $C_nS$ is disconnected with $k$ components if and only if 
the greatest common divisor of $\{n\} \cup S$ is $k$.  In particular, 
the disjoint union of two copies of $C_nS$ is a circulant graph 
$C_{2n}2S$ where $2S$ denotes the set $\{2s: s\in S\}$.

(C4) The connected graph $C_nS$ is $k$-partite if and only if $n$ 
is divisible by $k$ but no element of $S$ is divisible by $k$; in 
particular, a connected bipartite circulant graph has $n$ even and 
all elements of $S$ odd.

Similarly, we will state some well-known facts about the tensor product; 
most of these facts can be found in \cite{PMW, RobinThesis, MyThesis, 
Robin1, Robin2}.

(T1) The automorphism group of $G \otimes H$ always includes the direct 
product $\Aut(G) \times \Aut(H)$; the automorphism group of $G \otimes 
G$ always contains the wreath product $\Aut(G) \wr S_2$, with equality 
when $G$ is complete (\cite{Robin1, Robin2} and \cite{Handbook}, p.
1466).

(T2) The complement ${\cal C} (K_m\otimes K_n)$ is isomorphic to the 
Cartesian product $K_m \times K_n$.

(T3) A tensor product of connected graphs is disconnected if and only 
if both of the graphs are bipartite; then the product has exactly two 
components (\cite{PMW}).

(T4) A tensor product of two connected bipartite graphs has isomorphic 
components whenever one of the two graphs has an automorphism that maps 
one of its parts to the other.  (This condition is known to be 
sufficient and is conjectured to be necessary; the interested reader 
should see \cite{Jha1} and \cite{RobinThesis}.)

(T5) A tensor product of graphs is bipartite if and only if one of the 
graphs is bipartite.

We first state a basic result originally shown for all the usual graph 
products by Broere and Hattingh (\cite{B&H}).
%

\begin{theorem} If $G$ and $H$ are circulant graphs, and 
$(V(G),V(H)) = 1$, then $G \otimes H$ is a circulant graph.
\end{theorem}

\vskip 0.1 true in

This result establishes a condition that's sufficient but not necessary 
in general; indeed, the example $K_2 \otimes K_2 = C_4\{2\}$ belies the 
converse.  However, we may show a partial converse.  
%

\begin{theorem} Let $m$ and $n$ satisfy $(m,n)>1$ but at least one of 
$m$ and $n$ different from 2.  Then $K_m \otimes K_n$ is not circulant.
\end{theorem}

Again, Broere and Hattingh (\cite{B&H}) have shown something equivalent 
to this, 
but only for the case in which $m$ and $n$ are distinct even numbers. 
The proof of this theorem will be split into two parts; first we take 
care of the case where $m \neq n$.
%

\begin{lemma} If $m \neq n$ and $(m,n)>1$ then $K_m \otimes K_n$ is not 
circulant.
\end{lemma}

\begin{proof} From \cite{Robin1} 
and \cite{Robin2} we know that $\Aut(K_m \otimes 
K_n)$ is just the Cartesian product $S_m \times S_n$.  By (C1) we know 
that it suffices to show that no element of this group has an order 
equal to $mn$.  If $(\sigma, \tau)$ is an element of $S_m \otimes S_n$, 
it cannot be transitive on the vertices of $K_m \otimes K_n$ unless 
$\sigma$ is transitive on the vertices of $K_m$ and $\tau$ is 
transitive on the vertices of $K_n$.  Thus the order of $\sigma$ is 
$m$ and the order of $\tau$ is $n$.  But the order of $(\sigma, \tau)$ 
is just $\lcm(m,n)$ and because $m$ and $n$ are not relatively prime 
this is strictly less than $mn$.
\end{proof}

Now it merely suffices to prove the following theorem.
%

\begin{theorem} If $n \not= 2$, then $K_n \otimes K_n$ is not
circulant.
\end{theorem}

\begin{proof} From (T1) the automorphism group 
of $K_n \otimes K_n$ is the wreath product $S_n \wr S_2$.
Recall that the \emph{wreath product} $S_n \wr S_2 $ 
is defined to be the set of all objects of the form 
$(\sigma; \tau_1, \tau_2)$ where $\sigma \in S_2$ and $\tau_i \in S_n$.  
The action of this object on a vertex $(i_1,i_2)$ of $K_n \otimes K_n$ 
is defined to be $(\tau_{\sigma(1)}(i_{\sigma(1)}), 
\tau_{\sigma(2)}(i_{\sigma(2)}))$.  

In order for $K_n \otimes K_n$ to be circulant there must be an 
automorphism $\alpha$ that is transitive on the vertices of 
$K_n \otimes K_n$.
If $\alpha$ is to be transitive on $K_n \otimes K_n$ then 
the order of $\alpha$ must be $n^2$; by our earlier argument 
for the case $m \neq n$, that implies that $\alpha$ cannot be 
a member of the cross product, i.e., $\alpha = (\sigma; \tau_1, \tau_2)$ 
where $\sigma$ is not the identity element of $S_2$.  Since $\sigma$
has order 2, the order of $\alpha$ must be divisible by 2; hence
$n$ must be even.  Observe that
$\alpha = (\sigma; \tau_1, \tau_2)$ has the property that 
$\alpha^2(i_1,i_2) = (\tau_2(\tau_1(i_1)), \tau_1(\tau_2(i_2))) $ so 
that $\alpha^2 = (e; \tau_2 \tau_1, \tau_1 \tau_2)$ and is a member of 
the ordinary cross product of groups.  

Write $\tau_2 \tau_1 = \gamma$ and $\tau_1 \tau_2 = \delta$; and let 
$\gamma = \gamma_1 \gamma_2 \dots$ and $\delta = \delta_1 \delta_2
\dots$
be the disjoint cycle decompositions of $\gamma$ and of $\delta$, 
respectively.  We explicitly include any fixed points as cycles of
length 1.
%
%
%
%
Now suppose that 
$\gamma$ contains a cycle $\gamma_1$ of length $k < n/2$. 
Without 
loss of 
generality let us suppose that the cycle $\gamma_1$ contains the element
1.  
Consider the orbit of the vertex $(1, 1)$ under the action of
$\alpha^2$; 
it must contain at least $n^2/2$ elements.  
But the orbit
of 1 under the action of $\gamma$ is $\{i_1,   i_2, \dots, i_k\}$ where 
$k < n/2$ by choice of $\gamma_1$, so every vertex of the orbit
under 
$\alpha^2$ of $(1,1)$ must have as first coordinate one of the $i$'s. 
 There are at most $kn < (n/2)n = n^2/2$ such vertices, contradicting
the 
transitivity of $\alpha$. The same argument works to show the result for 
$\delta$, using second coordinates.
So we now know that $\gamma$ and $\delta$ (which are conjugate
permutations) are either both single $n$-cycles or both a product of 
exactly two $n/2$-cycles.  
Hence the orbit of $1$ under $\gamma$ and the orbit of $1$ under $\delta$ 
have the same size---either $n$ or $n/2$.
%
And since the size of the orbit of $(1,1)$ under 
the action of $\alpha^2$ is the least common multiple of the sizes of the
orbits of $1$ under $\gamma$ and $\delta$, 
this implies that 
the orbit of $(1, 1)$ under the action of $\alpha^2$ is at most $n$, which
contradicts the assumption of transitivity of $\alpha$ on $K_n \otimes
K_n$.

Since $\alpha$ cannot be transitive, it follows that 
$K_n \otimes K_n$ is not circulant for any value of $n > 2$. 
\end{proof}

The tensor product $K_m \otimes K_n$, as observed earlier, is the 
complement of the Cartesian product $K_m \times K_n$.  The graph $K_2 
\otimes K_2$ is the circulant graph $C_4\{2\}$.  
When $m$ and $n$ are relatively prime, 
$K_m \otimes K_n$ is easily seen to be isomorphic to 
$C_{mn}\{i: m \not|\ i, n \not|\ i\}$; that is, $S$ is the set of
integers 
between 1 and $\lfloor mn/2 \rfloor - 1$ that are not multiples of $m$ 
or $n$.  
%

\begin{corollary} The Cartesian product of two complete graphs $K_m 
\times K_n$ is circulant if and only if either $(m,n) = 1$ or 
$ m = n = 2$. 
\end{corollary}

\begin{proof} Immediate from (C2) and (T2) and Theorem~2. \end{proof}

It now can be shown that when $m$ and $n$ are relatively prime, $K_m 
\times K_n$ is isomorphic to $C_{mn}\{i: m | i \mbox{ or } n | i\}$.  
Additional results concerning the Cartesian and other products will be 
shown in a forthcoming paper by Sanders (``Products of circulant 
graphs are metacirculant'').

Next we introduce a result that will be useful in describing other 
tensor products as circulant graphs, as it invariably arises in 
components of disconnected products of connected bipartite graphs.  
Denote by $K_n^*$ the complete graph on $n$ vertices with a loop at 
each vertex; this is the graph with adjacency matrix consisting of all 
1's.
%

\begin{lemma}Let $G$ be a circulant graph on $m$ vertices.  Then
$K_n^* \otimes G$ is circulant.
\end{lemma}

\begin{proof}
Let the vertices of $G$ be labeled $0$,
$1$, \dots, $m-1$ in such a way that the permutation $(0, 1, \dots, m-1)$
belongs to $\Aut(G)$.  Note that two vertices $i$ and $j$ of $G$ are 
adjacent
if and only if $0$ is adjacent to both $j-i$ and $i-j$ where the 
subtraction is modulo $m$.  It immediately follows that $(i, j) \sim 
(k, \ell)$ in $K_n^* \otimes G$ if and only if $j \sim \ell$ in $G$.  

Define a map $\alpha \!:\! V(K_n^* \otimes G) \to V(K_n^* \otimes G)$\
as follows:
\[ \alpha(i, j) = \left\{
\matrix{
(i, {j+1}) \ \hbox{\rm if}\ j \not= m-1 \cr
(i+1, 0) \ \hbox{\rm if}\ j = m-1}
\right.
\]  
We claim that $\alpha$ is an automorphism of $K_n^* \otimes G$ and
$\alpha$ is a single cycle of order $nm$.  Once those facts are 
established, it clearly follows that $K_n^* \otimes G$ is a circulant 
graph with $nm$ vertices.

To see that $\alpha$ is a $nm$-cycle, consider the cycle of $\alpha$
that
contains the vertex $(0, 0)$: 
\begin{eqnarray*} 
& &\bigr( (0, 0), (0, 1), (0, 2), \dots, (0, {m-1}), (1,0), (1,1),
\dots,(1, {m-1}), 
\\ & & (2, 0), (2, 1), \dots, (n-1, 0), (n-1, 1), \dots, (n-1, {m-1})
\bigl)
\end{eqnarray*}

A better way to see this is to arrange the vertices of $K_n^* \otimes G$
into $n-1$ circles (layers) each containing $m$ vertices.  Think of the 
$n-1$ vertices that are lined up on top of each other as a ``slice.''  
The vertex $(i, j)$ belongs to the $i^{th}$ layer and the $j^{th}$
slice.  Then $\alpha$ simply moves a vertex in any slice other than
slice $m-1$ to the vertex on the same layer, but in the next slice, and
$\alpha$ takes a vertex in the last slice $(i, {m-1})$ to the
``starting'' 
vertex in the next layer---the vertex in slice 0  and layer $i+1$.  When
considered this way, it's clear that $\alpha$ is a $mn$-cycle of
vertices.

So now we must show that $\alpha$ is a graph automorphism of $K_n^*
\otimes G$.  
Suppose that $(i, j) \sim (k, \ell)$.
There are three cases to consider.  If neither $j$ nor $\ell$ equals 
$m-1$, then $\alpha(i,j) = (i, j+1)$ and $\alpha(k, \ell) = (k, \ell +
1)$.    
It is not hard to see that $\alpha$ will preserve the edge in this case
because $G$ 
is circulant.

If both $j$ and $\ell$ equal $m-1$, then $\alpha(i,m-1) = (i+1, 0)$, and 
$\alpha(k, m-1) = (k+1, 0)$.  Now, since $(i,m-1)$ is adjacent to $(k,
m-1)$, 
by the definition of the tensor product $G$ must have a loop at vertex
$m-1$.  
Since $G$ is circulant, $G$ must also have a loop at $0$, so $(i+1, 0) 
\sim (k+1, 0)$.


If one of $j$ or $\ell$ (say $\ell$) is equal to $m-1$, then $\alpha$ 
sends $(i,j)$ to $(i, j+1)$ and $(k, m-1)$ to $(k+1, 0)$.  Now, the 
completeness of $K_n^*$ implies that $i \sim k+1$; and since $j 
\sim m-1$ in $G$ and $G$ is circulant, we have $j+1 
\sim 0$ in $G$.  Therefore $(i, j+1) \sim (k+1, 0)$.

In all the cases, $\alpha$ preserves edges.  Thus $\alpha$ is a 
graph automorphism of $K_n^* \otimes G$.  And by (C3), $K_n^* \otimes G$ 
is circulant. 
\end{proof}
%

\begin{theorem} Let $G$ be a circulant bipartite graph.  Then $K_{n,n} 
\otimes G$ is circulant.
\end{theorem}

\begin{proof} It is straightforward to confirm that the graph $K_{n,n} 
\otimes G$ is made up of two disjoint copies of the graph $K_n^* 
\otimes G$.  (That this product is disconnected comes from (T3); the 
rest of this assertion is proven rigorously in \cite{RobinThesis}.)  
By Lemma~6, these components are circulant.  Now, two disjoint copies 
of a circulant graph form a circulant graph, by (C3), and we are done.
\end{proof}

The foregoing result may be extended somewhat.  Note, using (C4), that 
a circulant graph with an odd number of vertices cannot be bipartite.
%

\begin{theorem}
Let $G$ be a circulant graph with an odd number of vertices.  Then 
$K_{n,n} \otimes G$ is circulant.
\end{theorem}

\begin{proof}
The graph $K_{n,n}$ is isomorphic to $K_n^* \otimes K_2$, so 
$K_{n,n} \otimes G$ is isomorphic
to $(K_n^* \otimes K_2) \otimes G$ which (since the tensor product 
is associative) is isomorphic to  $K_n^* \otimes (K_2 \otimes G)$. 
Since $G$ has an odd number of vertices, Theorem~1 implies
that $K_2 \otimes G$ is circulant.  Lemma~6 implies that
$K_n^* \otimes (K_2 \otimes G)$ is circulant, and we're done.
\end{proof}

Clearly, if $G$ is any graph such that $K_2 \otimes G$ is circulant,
the above proof can be used to show $K_{n,n} \otimes G$ is circulant.

The next obvious question then becomes, ``When is $K_2 \otimes G$ 
circulant?''  Clearly, $K_2 \otimes G$ is circulant when $G$ is 
circulant and has an odd number of vertices
or when $G$ is circulant and bipartite.  
But the case where $G$ is circulant, non-bipartite, and has an even
number of vertices is more complex.  
For example, $K_2 \otimes K_4$ is equal to $Q_3$, the three dimensional
cube,
and hence is not circulant.  But if $G = C_6\{1,2\}$, then 
$K_2 \otimes G$ turns out to be isomorphic to $K_2^* \otimes C_6$ which
is 
circulant by Lemma~6.    
Strangely, it is possible for $K_2 \otimes G$ to be
circulant when $G$ is not.  For example, let $G$ be a path of length $n$
where 
the two end vertices, but not the internal vertices, have
a loop.  Then $K_2 \otimes G \cong C_{2n}$,
and is therefore circulant.  

Now we consider the general question of when the
product of two bipartite, circulant graphs is circulant.  A 
\emph{switching automorphism} of a bipartite graph is an automorphism 
that sends vertices from one part into the other.
%

\begin{theorem}
Let $G$ be a bipartite, circulant graph.  Then $G$ has a switching 
automorphism of order~2 and may be written as $K_2 \otimes H$ for some 
$H$ that may contain loops.
\end{theorem}

\begin{proof}
Observe that every circulant graph on $2n$ vertices has within its 
automorphism group a dihedral group $D_{2n}$ of order $4n$.  We label 
the vertices so that the cycle $\alpha = (0, 1, 2, \dots, 2n-1)$ is 
an automorphism.  Now the partition sets are $E =\{0,2,\dots,2n-2\}$ 
and $O = \{1,3,\dots,2n-3\}$.  The dihedral group is generated 
by $\alpha$ and $\beta = (0,2n-1)(1,2n-2)\dots(n,n-1)$.  By definition, 
$\beta$ is a switching automorphism.  The conclusion that $G = K_2 
\otimes H$ for appropriate $H$ is proven in \cite{gpw}.
\end{proof}
%

\begin{theorem} Let $G$ and $H$ be two connected, circulant, bipartite 
graphs on $2n$ and $2m$ vertices respectively.  If $(n,m) = 1$, then 
$G \otimes H$ is circulant.
\end{theorem}

\begin{proof}
By the previous theorem, both $G$ and $H$
contain switching automorphisms.   By (T4),  
both the connected components of $G \otimes H$ are isomorphic.  Hence,
using (C3), $G \otimes H$ is circulant if and only if one of its 
components is circulant.  Each component of $G \otimes H$ has $2nm$ 
vertices in it.  Now let $\sigma$ be a cycle of order $2n$ in $\Aut(G)$ 
and $\tau$ be a cycle of order $2m$ in $\Aut(H)$.  Then the element
$(\sigma, \tau)$ is an element of order $2nm$ in $\Aut(G \otimes H)$.
Once we show that $(\sigma, \tau)$ acts as a transitive graph
automorphism
on one of the  connected components
of $G \otimes H$, we're done.

The vertices of $G$ and $H$ may be labeled so that the 
partition sets are $E_G = \{g_{2i}\}$, $O_G = \{g_{2i+1}\}$,
$E_H=\{h_{2i}\}$, and $O_H= \{h_{2i+1}\}$.  Without loss of 
generality, we can also assume that $\sigma = (g_0, g_1, \dots,
g_{2n-1})$
and $\tau = (h_0, h_1, \dots, h_{2n-1})$.  Note that the first
connected component (Component~1) of $G \otimes H$ has vertex set
$(E_G \times E_H) \cup (O_G \times O_H)$.  Since
$(\sigma, \tau)$ sends $(g_i, h_j)$ to $(g_{i+1}, h_{j+1})$, it is clear
that $(\sigma, \tau)$ sends any vertex in $E_G \times E_H$ to 
a vertex in $O_G \times O_H$ and vice versa.  Hence, $(\sigma, \tau)$
acts as a graph automorphism on Component~1.  Transitivity of
$(\sigma, \tau)$ follows since each of $\sigma$ and $\tau$ are
transitive on their respective graphs.  Thus Component~1 is
circulant and so is $G \otimes H$.
\end{proof}

Now we demonstrate a large family of tensor products of graphs that are 
not circulant.  It is worthwhile to notice that this theorem even 
applies when the product is a product of circulant graphs.
%

\begin{theorem} Let $G$ be connected, bipartite, and regular of 
odd degree.  Let $H$ be connected, regular of odd degree, and not 
bipartite.  Then $G \otimes H$ is not circulant. 
\end{theorem}

\begin{proof} Observe $G \otimes H$ is regular of odd degree.  Since 
both graphs are regular of odd degree, it must be that 
$|G \otimes H| = 4n$ for some integer $n$.  
By (T5), the product is bipartite.  Now if a circulant graph $C_{4n}S$ 
is regular of odd degree then $S$ must contain $2n$.  However, by 
(C4) such a graph cannot be bipartite and so is not isomorphic to 
$G \otimes H$.
\end{proof}

We conclude with a discussion of open problems and interesting areas 
for further research.  We have left almost untouched the question of 
circulant products of non-circulant graphs.  Can one characterize those 
non-circulant graphs $G$ for which $K_{n,n} \otimes G$ is circulant?  
Such graphs do exist, by the remark following Theorem~8; are the graphs
shown 
there the only examples?  Also, if $G$ is circulant, not bipartite, 
and has an even number of vertices, what conditions on $G$ guarantee 
that $K_{n,n} \otimes G$ is circulant?  

Another question that arises from the proof of Theorem~2 concerns 
automorphism groups; in that proof, we observed that the automorphism 
group of $K_n \otimes K_n$ is $S_n \wr S_2$, which group is abstractly 
isomorphic to the automorphism group of $2K_n$, the graph consisting 
of two disjoint copies of $K_n$.  Now, $2K_n$ is circulant on $2n$ 
vertices, and $K_n \otimes K_n$ has $n^2$ vertices.  Is it possible 
to say that whenever a graph $G$ has more vertices than a circulant 
graph $H$, and $G$ and $H$ have isomorphic automorphism groups, that 
$G$ cannot be circulant?  This would assert, roughly, that of all graphs 
``realizing'' a certain abstract group as an automorphism group, a 
circulant graph would be one with the fewest possible vertices.


\begin{thebibliography}{44}
\bibitem{B&M} J.A.\ Bondy and U.S.R.\ Murty, \textsl{Graph Theory with 
Applications}, North-Holland, New York, 1976.

\bibitem{B&H} I.\ Broere and J.H.\ Hattingh, Products of circulant
graphs, 
\textsl{Quaestiones Mathematicae} \textbf{13}(1990), pp. 191--216.


\bibitem{MyThesis} J.C.\ George, One-factorizations of tensor products 
of graphs, Ph.D. thesis, University of Illinois at Urbana-Champaign,
1991.

\bibitem{gpw} J.C.\ George, T.D.\ Porter, and W.D. Wallis,
Characterizing 
balanced bipartite graphs with part-switching automorphisms, 
to appear in \textsl{BICA}.

\bibitem{Handbook} R.L.\ Graham, M.\ Gr\"otschel, and L.\ Lov\'asz, 
\textsl{Handbook of Combinatorics} v. II, The MIT Press, Cambridge, 
Massachussetts, 1995.

\bibitem{H&W} Frank Harary and Gordon W.\ Wilcox, Boolean operations 
on graphs, \textsl{Math. Scand.}, \textbf{20} (1967), pp. 41--51.

\bibitem{Jha1} Pranava K.\ Jha, Sandi Klav\v{z}ar, and Bla\v{z} Zmazek, 
Isomorphic components of Kronecker product of bipartite graphs. 
\textsl{Discuss. Math. Graph Theory}, \textbf{17} (1997) no. 2, pp.
301--309.

\bibitem{RobinThesis} R.S.\ Sanders, Graphs on which dihedral,
quaternion, 
and abelian groups act vertex and/or edge transitively and applications 
to tensor products, Ph.D. thesis, University of Illinois at
Urbana-Champaign, 1990.

\bibitem{Robin1} R.S.\ Sanders and J.C.\ George, Basic results 
concerning the automorphism group of the tensor product of two graphs, 
\textsl{Utilitas Mathematica}, June, 1997 pp. 51--63.

\bibitem{Robin2} R.S.\ Sanders and J.C.\ George, Results concerning the 
automorphism group of the tensor product $G \otimes K_n$,
\textsl{JCMCC}, 
\textbf{24}(1997), pp. 119--127.

\bibitem{PMW} P.M.\ Weichsel, The Kronecker product of graphs, 
\textsl{Proc. Amer. Math Soc.} \textbf{13}(1962), 47--52.

\end{thebibliography}
\end{document}